\documentclass[conference]{IEEEtran}
\ifCLASSINFOpdf
\else
\fi
\usepackage[font=small,margin=10pt,labelfont={bf},labelsep={space}]{caption}
\usepackage{subfig}
\usepackage{wrapfig}
\usepackage{amsmath, amssymb, epsfig}
\usepackage[scaled]{helvet} % I like Helvetica for sf
\usepackage{fourier}
\usepackage{bm}
\usepackage{color}
\usepackage{fullpage} 	% Fullpage package

\newcommand{\hide}[1]{}

\newcommand{\ignore}[1]{}

\newcommand\remove[1]{}

\usepackage{palatino}

\usepackage{hyperref}

\usepackage[T1]{fontenc}
\usepackage[utf8]{inputenc}
\usepackage{tabularx,ragged2e,booktabs,caption}

\usepackage{algorithmicx}
\usepackage[ruled]{algorithm}
\usepackage{algpseudocode}

\begin{document}
\title{A Comparison of Clustering and Missing Data Methods for Health Sciences}

\author{
\IEEEauthorblockN{R. Zhao}
\IEEEauthorblockA{Institute of Mathematical Sciences \\
Claremont Graduate University \\
Claremont, CA, 91711 \\
Email: ran.zhao@cgu.edu }%\\
\and
\IEEEauthorblockN{D. Needell}
\IEEEauthorblockA{Department of Mathematics\\
and Computer Science \\
Claremont McKenna College \\
Claremont, CA, 91711 }
\and
\IEEEauthorblockN{C. Johansen and J. L. Grenard}
\IEEEauthorblockA{School of Community \\ and Global Health\\
Claremont Graduate University \\
Claremont, CA, 91711}
}


%Telephone: (800) 555--1212, Fax: (888) 555--1212}

\thanks{Correspondence Information: Ran Zhao, Institute of Mathematical Sciences at Claremont Graduate University, Claremont, CA 91711, tel: (909) 969-0553, Email:ran.zhao@cgu.edu}


\maketitle

\begin{abstract}
%\boldmath
In this paper, we compare and analyze clustering methods with missing data in health behavior research.  In particular, we propose and analyze the use of compressive sensing's matrix completion along with spectral clustering to cluster health related data.  The empirical tests and real data results show that these methods can outperform standard methods like LPA and FIML, in terms of lower misclassification rates in clustering and better matrix completion performance in missing data problems. According to our examination, a possible explanation of these improvements is that spectral clustering takes advantage of high data dimension and compressive sensing methods utilize the near-to-low-rank property of health data.
\end{abstract}


% \fi
\IEEEpeerreviewmaketitle

\section{Introduction}
\subsection{Clustering Analysis}
A vast array of literature has explored clustering techniques and missing data issues in both mathematics and public health research. Clustering refers to the separation of data into meaningful groups so that data within each group is similar.

The Latent Profile Analysis (LPA) method is a common approach in health behavior research to identify unobserved classes of participants and explain the pattern of responses~\cite{Vermunt02, Lubke05, Wang12, Marsh09}. Many current software packages use an iterative expectation maximization (EM) algorithm to estimate the parameters~\cite{Mplus}. The EM algorithm and other variants have both advantages and drawbacks for estimation of the LPA parameters. The algorithms are sensitive to the initial values of the parameters with the potential for local solutions, and the EM approach does not estimate standard errors. Model identification, the issue of whether there is sufficient information to estimate the parameters~\cite{Vermunt02}, and subjective model fit selection are also drawbacks to these approaches.

Spectral clustering (SC) is a geometric method that can identify relationships in the data (here we consider $n$ individuals each with $d$ variables) that are non-linear~\cite{lloyd1982least, shi2000normalized, wu2009top}. Here, one designs a similarity measure to form a \emph{Laplacian} matrix from the data.  A typical normalized \emph{Laplacian} matrix $\mathbf{L}\in\mathbb{R}^{n\times n}$ is defined by

\begin{equation} \label{eqn1}
L = D^{-1/2}(D-W)D^{-1/2},
\end{equation}

where $W$ is the symmetric weight matrix whose $(i,j)$th entry corresponds to the similarity between individuals $i$ and $j$, and the degree matrix $D$ has diagonal entries $D_ii = \sum_j W_{ij}$.  Spectral clustering computes the eigenvectors of this Laplacian which form a lower dimensional, linear separable representation of the dataset~\cite{Shi00}.

\begin{figure}[ht]
\includegraphics[width=3.5in]{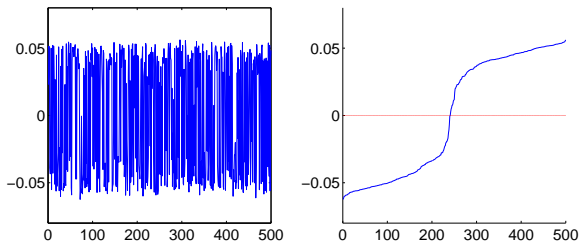}
\caption{Left: unsorted eigenvector.  Right: Sorted eigenvector.  The red horizontal line indicates the separation. \label{fig-1}}
\end{figure}

%These low dimensional vectors are quite easy to observe.
In the dataset from Section~\ref{sec:clu}, the unsorted and sorted eigenvector (from the second largest eigenvalue) entries are shown in Figure~\ref{fig-1}.  Since each entry of the eigenvector corresponds to an individual, we may use the values of these entries to separate the individuals into clusters.  In this case, the threshold to designate different clusters seems to be $y = 0$, as plotted in red. For more than two clusters, one can instead run $k$-means~\cite{lloyd1982least,wu2009top} on this data to identify the clusters.

% \dn{Need to fix labels and make bigger.}

\begin{figure}[ht]
\includegraphics[width=3in]{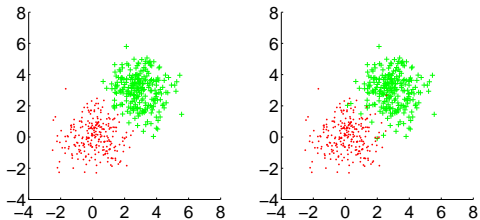}
\caption{Results from spectral clustering (on the left) actual real clusters (on the right). Red and green show different clusters. \label{fig0}}
\end{figure}

%According to the clusters results of the key eigenvector, we then distinguish clusters in real observations. In the same case, spectral clustering separate the observations into two groups very well.
In Figure~\ref{fig0}, we show the comparison of the spectral clustering results and that of the actual (randomly generated) clusters on right.

% \dn{Need to fix labels and make bigger.}

\subsection{Missing Data}
In many large scale applications, data is incomplete. For example, participants may be unable or unwilling to complete an ongoing survey, or participants may be randomly assigned different blocks of questions to increase the variety of constructs assessed.

\subsubsection{FIML}
Full Information Maximum Likelihood (FIML)~\cite{Graham12,Enders10} aims to maximize the likelihood of the data by auditioning combinations of parameter estimates~\cite{Ender01}. The procedure relies on assumptions such as normality, which when violated can result in biased parameters. There is also a risk of convergence to local maxima resulting in poor parameter estimates. %The current study applied the FIML approach to adjusting for missing data.

The FIML estimator implemented in common statistical software packages maximizes a likelihood function that is the sum of $n$ case-wise likelihood functions. A likelihood function is calculated for each observation or individual. The function measures the discrepancy between the current parameter estimates and the observed data for the $i$th case. The function is maximized assuming multivariate normality:

\begin{equation*}
\log L_i = K_i - \frac{1}{2}\log |\Sigma_i| - \frac{1}{2}(x_i-\mu_i)' \Sigma_i^{-1} (x_i-\mu_i)
\end{equation*}

The vector of complete data for case $i$ is represented by the term $x_i$, and the vector of estimated means for those variables that are observed for case $i$ is the term $\mu_i$. A constant, $K_i$ depends upon the number of complete data points for case $i$. Only those variables that are observed for case $i$ are used to calculate the determinant and inverse of $\Sigma_i$. The discrepancy function for the entire sample is calculated by summing over the $n$ case-wise functions:

\begin{equation*}
\log L(\mu,\Sigma) = \sum_{i=1}^{N}\log L_i.
\end{equation*}

It is assumed that missing values for X are conditionally dependent on other observed variables in the data. Probability values for the missing data are implied during the parameter estimation process by incorporating vectors of partially complete data in the individual-level likelihood functions. This is analogous to using multiple regression of X on other variables to generate predicted scores for the missing data. The FIML estimate does not impute missing values, however, but uses all available raw data to directly estimate parameters and standard errors for the model.

\subsubsection{Compressive Sensing}
Compressive sensing (CS) is a new and fast growing field in applied mathematics. The CS application \textit{matrix completion} demonstrates that a (nearly) low-rank matrix can be completed accurately and robustly from observation of only a few of its entries by solving a nuclear-norm minimization problem~\cite{DSPweb,CandeRT_Stable,Donoho06}. A typical format of this optimization problem is

\begin{eqnarray*}
\textrm{minimize} & ||X||_{*} \\
\textrm{subject to} & X_{ij} = M_{ij} \quad (i,j) \in \Omega
\end{eqnarray*}
where the nuclear norm $||X||_* = \sum_{k=1}^{n}\sigma_k(X)$, $M$ is the matrix we wish to recover, and $\Omega$ is the set of locations of observed matrix entries in $M$. This popular convex relaxation of the rank minimization problem is feasible and commonly used in matrix completion, since minimization of the rank of $X$ is NP-hard due to its combinatorial nature.

When the underlying matrix is low-rank, matrix completion completes the data matrix provably well from a small number of possibly noisy observations~\cite{candes2010matrix,nep:rank,recht2007guaranteed}. To quantify how exact the method recovers the matrix, we generate two $1000\times 300$ matrices with rank $2$ and $10$ and remove $20\%$, $40\%$, $60\%$, and $80\%$ of the data purposefully. The entries are random values that follow a standard normal distribution. The matrix completion results are presented in Table~\ref{table-1} and Table~\ref{table0}. We measure the recovery error between the actual matrix $X$ and the recovered matrix $\hat{X}$ by the Frobenius norm $\|X-\hat{X}\|_F$, the relative Frobenius norm $\|X-\hat{X}\|_F/\|X\|_F$, and the spectral norm $\|X-\hat{X}\|$.  As is evident and not surprising, the error increases slightly with more missing data, and the higher rank matrix has slightly higher recovery error.  %It is obvious that error matrix with lower rank has lower Frobenius norm, indicating better matrix completion. %The Frobenius norm of the error matrix can be interpreted as the ``distance'' between the recovered and the original matrix in their column space. We observe a very small ``distance'' or difference from the original matrix after completing by compressive sensing.

\begin{table}
\begin{minipage}{\linewidth}
\centering
% Table 1
\caption{Rank 2 Matrix Completion Results}
%In this table, matrix completion measurements are shown for the rank 2 matrix. Frobenius represents the Frobenius norm of the different between recovered matrix and the original one, which is also called error matrix. Relative Frob. represents the Frobenius norm of error matrix divided by that of the original matrix. Spectral at the top of last column represents the spectral norm of the error matrix.
\begin{tabular}{cccc} \toprule[1.5pt]
Rank 2 &	Frobenius &	Relative Frob. &	Spectral \\ \hline
missing 20\% &	0.0687 &	8.84E-05 &	0.0505 \\
missing 40\% &	0.0376 &	4.85E-05 &	0.0246 \\
missing 60\% &	0.0651 &	8.38E-05 &	0.0422 \\
missing 80\% &	0.0959 &	1.23E-04 &	0.0645 \\
\bottomrule[1.25pt]
\end{tabular}\par
\label{table-1}

\bigskip

% Table 2
\caption{Rank 10 Matrix Completion Results}
%In this table, matrix completion measurements are shown for the rank 10 matrix. All column names are the same as that of Table~\ref{table-1}.
\begin{tabular}{cccc} \toprule[1.5pt]
Rank 10 &	Frobenius &	Relative Frob. &	Spectral \\ \hline
missing 20\% &	0.0896 &	5.37E-05 &	0.0297 \\
missing 40\% &	0.145 &	8.69E-05 &	0.0512 \\
missing 60\% &	0.186 &	1.12E-04 &	0.0643 \\
missing 80\% &	0.350 &	2.10E-04 &	0.1890 \\
\bottomrule[1.25pt]
\end{tabular}\par
%\bigskip
\label{table0}
\end{minipage}
\end{table}

In public health data, especially the data from surveys or investigations, one expects the data to be low-rank or approximately low-rank for certain variables, because there are a small number of underlying factors that influence specific human opinions and behaviors.

In this paper, we empirically investigate the use of matrix completion with spectral clustering to cluster incomplete data, and compare to standard FIML and LPA methods. From these studies, we find that the combination of compressive sensing and spectral clustering methods can offer better performance than standard methods currently used in health data research.  %Spectral clustering and compressive sensing outperform LPA and FIML respectively, by the empirical results from simulated data sets.
%Our research seems to support the claim that spectral clustering and compressive sensing can together generate better outcomes in health behavior research than their counterparts, LPA and FIML, which are currently being widely and commonly used in public health data oriented research.

\section{Empirical Results}
\subsection{Experiments on Clustering Analysis}  \label{sec:clu}
We first use simulated data to compare the clustering performance of spectral clustering and LPA. First, we generate two-dimensional points whose $x$ and $y$ values follow a normal distribution with mean $0$ (for one cluster) or $a>0$ (for the other cluster) and variance $1$.  As $a$ increases, we expect clustering to be more successful, because the difference between clusters is more obvious.  %The following chart shows the proportion of correctly clustered points for spectral clustering and LPA with $a$  = $1$, $2$, $3$, and $5$.
We define the correct classification rate (CCR) as the ratio of the number of correctly clustered points over the total number of points in each trial.

\begin{figure}[ht]
\includegraphics[width=3in]{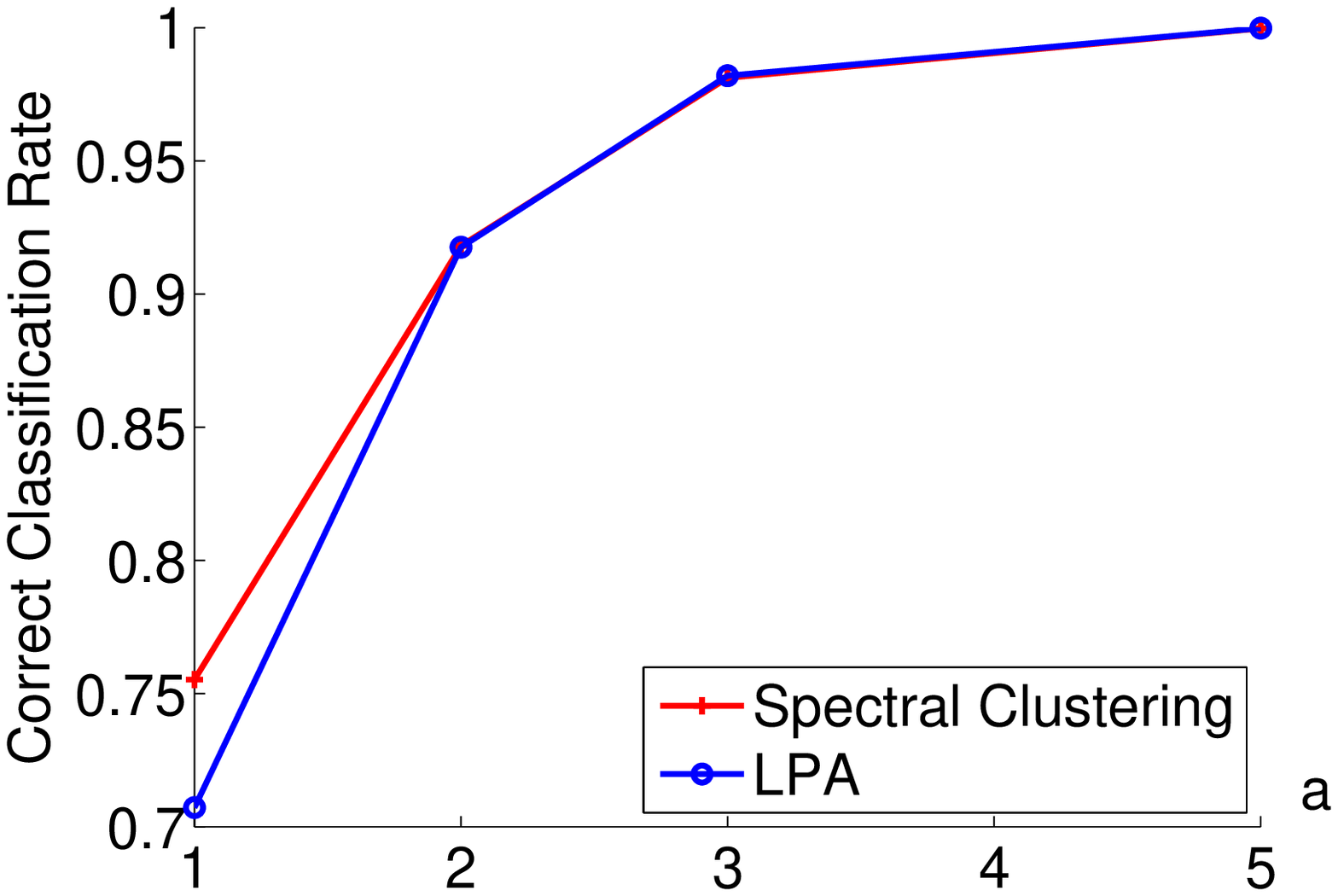}
\caption{Clustering results of spectral clustering and LPA methods for equally sized clusters.  \label{fig1-1}}
% \dn{Need to fix labels and make bigger.}
\bigskip

\includegraphics[width=3in]{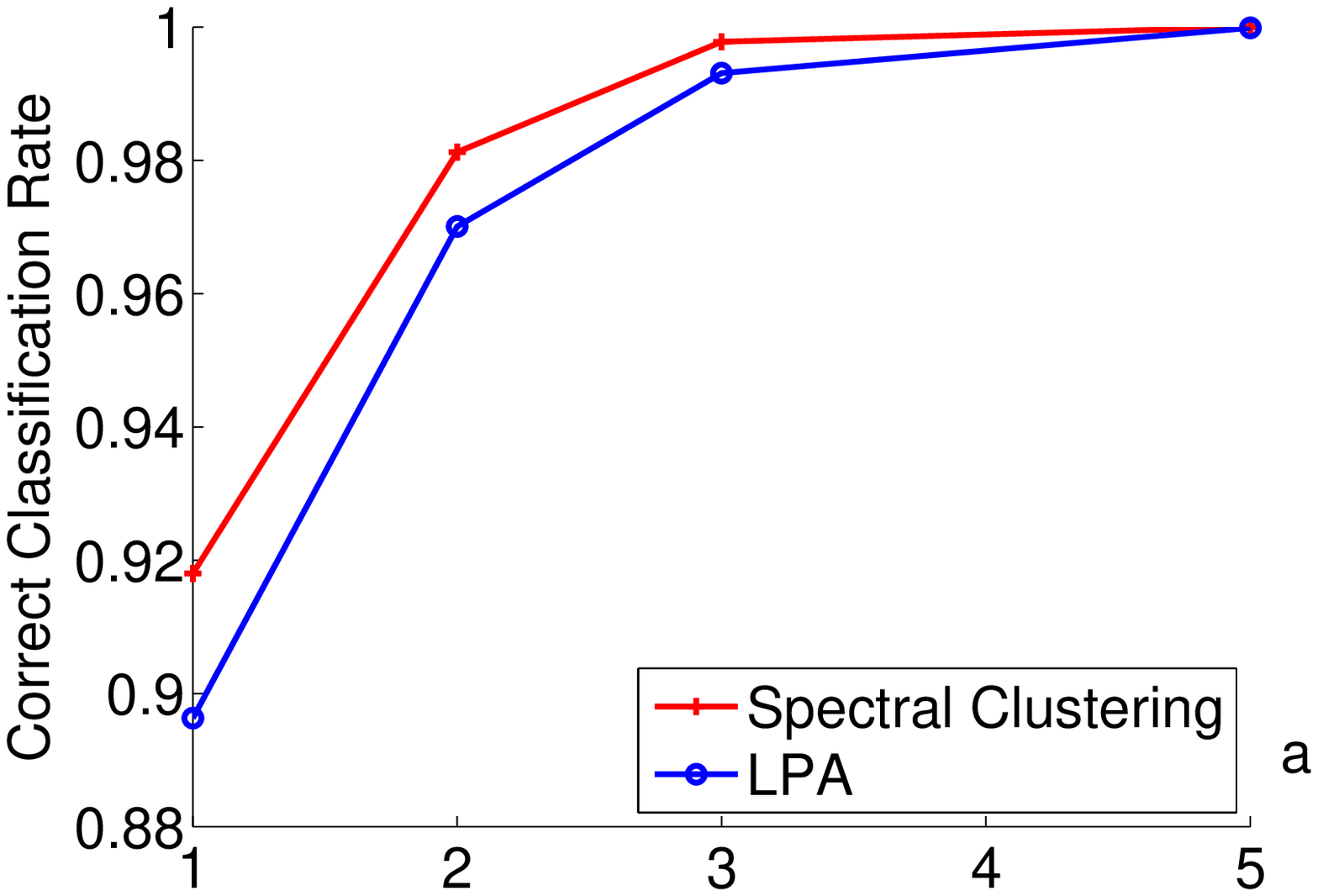}
\caption{Clustering results of spectral clustering and LPA methods for unequally sized clusters.  \label{fig1-2}}
\end{figure}

We simulate 40 different data sets for each value of $a$, and use these 40 trials to compute the rates of each method. %The mean of rate for each value of $a$ and each method is plotted in Figure~\ref{fig1-1} and Figure~\ref{fig1-2}.
In Figure~\ref{fig1-1}, we show the mean CCR for datasets that contain two equally sized clusters, with approximately $250$ observations in each cluster. Figure~\ref{fig1-2}, however, illustrates the CCR from datasets that contain unequal sizes of clusters, where one has approximately $25$ ($5\%$) observations and the other has approximately $475$ ($95\%$) observations. This experiment aims to test how well spectral clustering and LPA classify observations in situations with different clustering complexity and relative size of clusters.
 %\notate{Modified here}

%\dn{Need to fix labels and make bigger.}

\begin{table}
\begin{minipage}{\linewidth}
\centering
\caption{Correct Classification Rate Results}
In this table, $a$ represents the centroid $(a,a)$ other than $(0,0)$. CCR is the abbreviation of correct classification rate, and N represents the number of observations. In this test we have sample size of 500 and approximately 250 obs. for each cluster. The rest of items are summary statistics of correct classification rate estimated from 40 trials.
\begin{tabular}{cccccccc} \toprule[1.5pt]
$a$	&	CCR	&	N	&	Mean	&	S.D. &	Min	&	Med.	&	Max	\\	\hline
1	&	SC	&	500	&	0.755	&	0.0198	&	0.692	&	0.756	&	 0.812	\\	
	&	LPA	&	500	&	0.707 &	0.0550 &	0.508 &	0.072 &	0.796	\\	\hline
2	&	SC	&	500	&	0.918	&	0.0132	&	0.868	&	0.918	 &	0.958	\\	
	&	LPA	&	500	&	0.918 &	0.0130 &	0.882 &	0.916 &	0.962	\\	\hline
3	&	SC	&	500	&	0.981	& 0.0065	&	0.958	&	0.982	 &	0.998	\\		
	&	LPA	&	500	&	0.982 &	0.0060 &	0.960 &	0.982 &	0.998	\\	\hline
5	&	SC	&	500	&	1.00	& 0.0008	&	0.994	&	1.00	&	1	 \\	
	&	LPA	&	500	&	1.00 &	0.0010 &	0.996 &	1.00 &	1.00	\\	
\bottomrule[1.25pt]
\end{tabular}\par
% \bigskip
\label{table1}
\end{minipage}
\end{table}

We observe that both methods have increasing CCR for larger $a$, and that spectral clustering has a higher correction rate when the distance between centroids is small for equally sized clusters and for unequally sized clusters. Overall, spectral clustering seems to offer improvements over LPA in this setting. A detailed summary of results for equal cluster sizes is shown in Table~\ref{table1}. When $a\geq 2$, where the (spatial) distance between centroids is greater than $2.828$, the CCR of both methods exceed $90\%$, but in most cases spectral clustering has lower standard deviation of estimation.

\subsection{Experiments on Missing Data}\label{sec:missdata}
Missing data is an important part of real-world health data analysis.
There is an obvious difference in how FIML and compressive sensing handle this problem.  FIML does not recover missing data or interpolate unknown values to incomplete data entries, while matrix completion does precisely this. Given this difference, it is hard to compare the performance of FIML and compressive sensing directly, so we instead use their results to cluster.  We compare the correction rates of the following three methods, a) FIML combined with LPA, b) matrix completion followed by LPA, and c) matrix completion followed by spectral clustering.

By comparing methods a) and b), we can compare the relative performance of FIML and matrix completion, because the clustering methods are the same (LPA). Therefore the correct classification rates will reflect how well these two methods handle missing data.  The results of b) and c) will strengthen the conclusion from part a), where the two different clustering methods are compared.

We generate $40$ data sets of $1000$ (the number of individuals) by $100$ (the number of variables) matrices.  One can imagine the $500\times 100$ top half of the matrix corresponding to one cluster, and the bottom half to the other.  The bottom half has standard normally distributed entries, whereas the top half has normally distributed entries with variance $1$ but with varying means; the first ten columns have mean $0.1$, the next $10$ have mean $0.2$, and so on, so that the last ten have mean $1.0$ (we do this to introduce more variety within the cluster).  We remove entries from the matrix uniformly at random to create missing data. %The correction rates of the above three methods with randomly removed data are shown below.

\begin{figure}[ht]
\includegraphics[width=3in]{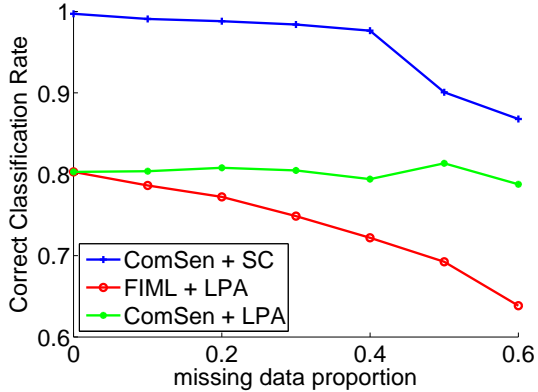}
\caption{Mean correct classification rate as a function of the proportion of missing data for method a), method b), and method c), (red, green, and blue, respectively). The description of the methods is elaborated in~\ref{sec:missdata}.
%The horizontal axis represents the proportion of missing data, the vertical represents the mean correction rate. The red line represent the results from method a), which is FIML followed by LPA; the green one illustrates results from method b), which is compressive sensing followed by LPA; the results  form method c), which is compressive sensing followed by spectral clustering, are plotted in blue. \notate{Does this make it clearer?}
\label{fig2}}
\end{figure}

In Figure~\ref{fig2}, we observe that the correct classification rate is higher for method b) than for method a). This indicates that matrix completion may handle missing data better than FIML in this setting.  Additionally, method c) is better than both, suggesting that matrix completion coupled with spectral clustering may offer even better performance.  There is of course a computational tradeoff, but for methods where accuracy is the priority, matrix completion coupled with spectral clustering may offer improved performance over standard methods.  %\notate{The interpretation here is correct!} %Given the same recovered data matrix, spectral clustering generates higher correction rates than LPA on average.

\subsection{Application to Real Public Health Data}
We next compare the above methods on real public health data.  The data is obtained from the teen California Health Interview Survey (CHIS) from 2009. CHIS is one of the largest surveys in the nation and is conducted and maintained by the UCLA Center for Health Policy Research and its collaborators. %CHIS has collected health information from Californians every two years since 2001 for all age groups.
CHIS obtains data via phone interviews on extensive health related items such as health status, health conditions, health-related behaviors, health insurance coverage, access to health care services, and other health and health related issues~\cite{CHIS}.

One major difficulty of analyzing the clustering techniques on real data is that there is not an obvious ground truth to which to compare. %Also, researchers will select different groups of variables in their researches according to the topics.
% So different clusters and measurements of clusters vary. For example, one will select total different variables when judging adiposity and mental distress. In mental distress research, the quantitative definition of distress is sometimes ambiguous.
%Typically, researchers who use this data in the health sciences identify particular variables of interest and ignore others.
%To solve these problems, we use as much information of the CHIS data.
To over come this, we first eliminate irrelevant variables such as individual's serial number and zip code.  This yields a data matrix with 3379 individuals and 144 variables. Then we apply both spectral clustering and LPA, and identify those individuals who were clustered in the same way by both methods. This left 2836 individuals as the ``consistent population'', which is 83.93\% of the original data. Next, we sample $1000$ individuals (without replacement) among this consistent population, and repeat this process $40$ times. In each of these $40$ trials, we randomly remove $10\%$, $30\%$, and $50\%$ of the entries to mimic missing data. Finally, we apply matrix completion/spectral clustering and FIML/LPA and compute the mean CCR using the consistent clusters as ground truth for each approach.

\begin{figure}[ht]
\includegraphics[width=3in]{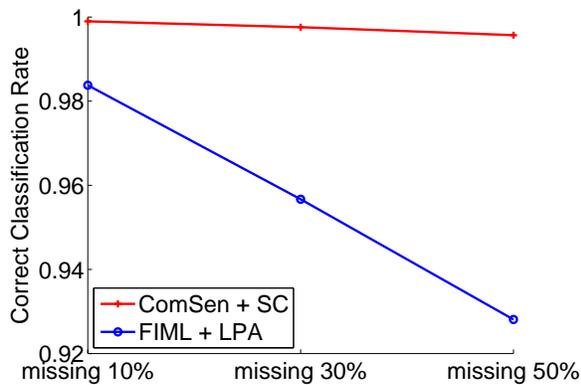}
\caption{Missing data completion and clustering analysis on CHIS data. %The horizontal line shows the proportion of missing data, and in our experiment they are 10\%, 30\% and 50\% respectively. The vertical line is final correct classification rates. The red line represents the mean of 40 trials by compressive sensing and spectral clustering combination. Whereas the blue line indicates the averaged results from FIML and LPA combination.
\label{fig3}}
\end{figure}

The averaged rates are illustrated in Figure~\ref{fig3}. Though as expected the correct classification rates decrease monotonically, the CCR for compressive sensing/spectral clustering seems to decay at a much slower speed than that of FIML/LPA. Regardless, these two approaches overall generate quite reliable outcomes, even when the proportion of missing data reaches as large as $50\%$.

\section{Summary}
Using two groups of simulated data, we observe that spectral clustering may be preferable to LPA, and that compressive sensing methods may have an advantage over FIML in giving the recovered data matrix explicitly, taking advantage of nearly low-rank data.

The contribution of this paper is to bring two methods from applied mathematics into health behavior research, and verify their advantages over traditionally used methods.  Our future research direction is to further compare the performance of these methods on real health (and other types of) data, and aim to identify in what settings each type of method is preferred.  This identification can aid in the design of health data surveys allowing for intentional missing data, thereby reducing participant burden and cost.

\section*{Acknowledgement}
Support for this research was provided by a BLAIS Challenge Award (221-2170045) received from Claremont Graduate University, Claremont, CA.



%\bibliographystyle{IEEEtran}

\bibliographystyle{plain}
\bibliography{bib}

\begin{thebibliography}{10}

\bibitem{DSPweb}
Compressed sensing webpage.
\newblock http://dsp.rice.edu/cs.

\bibitem{CHIS}
California health interview survey. chis 2009 adolescent public use file.
  release 1.
\newblock Technical report, UCLA Center for Health Policy Research, Los
  Angeles, CA, 2009.

\bibitem{CandeRT_Stable}
E.~Cand\`{e}s, J.~Romberg, and T.~Tao.
\newblock Stable signal recovery from incomplete and inaccurate measurements.
\newblock {\em Communications on Pure and Applied Mathematics},
  59(8):1207--1223, 2006.

\bibitem{candes2010matrix}
E.~J. Candes and Y.~Plan.
\newblock {Matrix completion with noise}.
\newblock {\em Proceedings of the IEEE}, 98(6):925--936, 2010.

\bibitem{Wang12}
S.J.H. Biddle W.C.~Liu C.K.J.~Wang and B.S.C. Lim.
\newblock A latent profile analysis of sedentary and physical activity paterns.
\newblock {\em Journal of Public Health}, 20:367--373, 2012.

\bibitem{Donoho06}
D.~Donoho.
\newblock Compressed sensing.
\newblock {\em IEEE Transactions on Information Theory}, 52(4):1289--1306,
  2006.
\newblock Published online.

\bibitem{nep:rank}
Y.~C. Eldar, D.~Needell, and Y.~Plan.
\newblock Uniqueness conditions for low-rank matrix recovery.
\newblock {\em Applied Compututational Harmonic Analysis}, 31(1):59--73, 2011.

\bibitem{Enders10}
C.~K. Enders.
\newblock {\em Applied missing data analysis}.
\newblock New York: Guilford Press, 2010.

\bibitem{Ender01}
C.K. Enders.
\newblock The performance of the full information maximum likelihood estimator
  in multiple regression models with missing data.
\newblock {\em Educational and Psychological Measurement}, 61(5):713--740,
  2001.

\bibitem{Graham12}
J.~W. Graham.
\newblock {\em Missing data analysis and design}.
\newblock New York: Springer, 2012.

\bibitem{Vermunt02}
J.~A. Hagenaars and A.~L. McCutcheon.
\newblock {\em Applied latent class analysis}.
\newblock Cambridge University Press, New York, 2002.

\bibitem{Marsh09}
U.~Tautwein H.W.~Marsh, L.~Ludtke and A.J.S. Morin.
\newblock Classical latent profile analysis of academic self-concept
  dimensions: Synergy of person- and variable-centered approaches to
  theoretical models of self-concept.
\newblock {\em Structural Equation Modeling}, 16(2):191--225, 2009.

\bibitem{lloyd1982least}
Stuart Lloyd.
\newblock Least squares quantization in pcm.
\newblock {\em Information Theory, IEEE Transactions on}, 28(2):129--137, 1982.

\bibitem{Lubke05}
G.H. Lubke and B.~Muthen.
\newblock Investigating population heterogeneity with factor mixture models.
\newblock {\em Psychological Methods}, 10(1):21--39, 2005.

\bibitem{Mplus}
L.K. Muthen and B.O. Muthen.
\newblock {\em Mplus User's Guide}.
\newblock Los Angeles, CA: Muthen \& Muthen, 1998-2012.

\bibitem{recht2007guaranteed}
B.~Recht, M.~Fazel, and P.A. Parrilo.
\newblock {Guaranteed minimum-rank solutions of linear matrix equations via
  nuclear norm minimization}.
\newblock {\em SIAM Review}, 52(3):471--501, 2010.

\bibitem{Shi00}
J.~Shi and J.Malik.
\newblock Normalized cuts and image segmentation.
\newblock {\em Pattern Analysis andMachine Intelligence, IEEE Transactions on},
  22(8):888--905, 2000.

\bibitem{shi2000normalized}
J.~Shi and J.~Malik.
\newblock Normalized cuts and image segmentation.
\newblock {\em IEEE Transactions on Pattern Analysis and Machine Intelligence},
  22(8):888--905, 2000.

\bibitem{wu2009top}
Xindong Wu and Vipin Kumar.
\newblock {\em The top ten algorithms in data mining}.
\newblock CRC Press, 2009.

\end{thebibliography}

%\begin{thebibliography}{1}
%\end{thebibliography}

\end{document}